\documentclass{article}

\newcommand\ov[1]{\overline{#1}}
\newcommand\ds[1]{\displaystyle{#1}}
\newcommand\csABt{\cos\left(\sqrt{A^2+B^2}\,t\right)}
\newcommand\ptl{\partial}
\newcommand\alp{\alpha}
\newcommand\lam{\lambda}
\newcommand\del{\delta}
\newcommand\Del{\Delta}
\newcommand\Gam{\Gamma}
\newcommand\ome{\omega}
\newcommand\Ome{\Omega}
\newcommand\sig{\sigma}
\newcommand\epsi{\varepsilon}
\newcommand\Ree{{\rm Re}$\,$}
\newcommand\Imm{{\rm Im}$\,$}
\newcommand\pll{\parallel}
\newcommand\rar{\rightarrow}
\newcommand\ppt{{\partial\over\partial t}}
\newcommand\izi{{\int_0^{\infty}}}
\newcommand\iii{{\int_{-\infty}^{\infty}}}
\newcommand\ltptr{{\left({1\over t}\ppt\right)}}
\newcommand\smboc{\quad\hbox{\hbox{\vrule width 4pt height 6pt depth 
  0pt}}} 
\def\eqalign#1{\null\,\vcenter{\openup\jot
  \ialign{\strut\hfil$\displaystyle{##}$&$\displaystyle{{}##}$\hfil
  \crcr#1\crcr}}\,}

\begin{document}
\title{The Method of Ascent and $\cos\left(\sqrt{A^2+B^2}\,
\right)$\footnote{This research was partially supported by a Minerva
foundation grant.}}
\author{Yakar Kannai\\ \\
Department of Mathematics\\
The Weizmann Institute of Science\\
Rehovot 76100, Israel}
\date{}
\maketitle
\section{Introduction}
The well-known Hadamard method of descent \cite{F,H}\ consists of 
deriving solutions of wave equations in $n$ variables from (known) 
solutions in more than $n$ variables. To quote Hadamard, (the method) 
``consists in noticing that he who can do more can do less''. It is the 
purpose of this paper to exhibit cases in which he who can do less can
do more --- that sometimes one can build the propagator $\csABt$ out of 
$\cos(At)$ and $\cos(Bt)$.\par
In Section~2 we consider the commutative case and represent\linebreak
$\cos\left(\sqrt{A_1^2+\cdots+A_n^2}\,t\right)$ in terms of an integral
involving $\cos(A_it_i)$ (Theorem~2.1). We distinguish between the 
cases $n$ even and $n$ odd. While it is true that once you see the 
formula ((2.1) or (2.2)), you may verify it directly (and this is the 
method used for the case of $n$ even), we prefer to use for the case $n$ 
odd a method which shows how the trigonometric-hyperbolic case is 
derived from the elementary exponential (parabolic) case $\exp\left(-
\rho(A_1^2+\cdots+A_n^2)\right)$ through ``inversion'' of a 
transmutation formula \cite{K}. The case $n$ even could also be reduced 
to the case $n$ odd via a descent argument. Note that (2.1) and (2.2) 
appear to be new even for the case of numbers (scalar operators). 
The cases $n=2,3$ are rather simple:
$$\csABt={1\over 2\pi}\ppt\bigg[t \int_{\ome_1^2+\ome_2^2} {\cos(t\ome_1
A)\cos(t\ome_2 B)\over\sqrt{1-(\ome_1^2+\ome_2^2)}}d\ome_1 d\ome_2
\bigg]\eqno (1.1)$$ 
$$\eqalign{&\cos\left(\sqrt{A^2+B^2+C^2}\,t\right)\cr
\noalign{\vskip 5pt}
&\qquad ={1\over 4\pi}\ppt\ltptr\bigg[t \int_{S^2}\cos(t\ome_1 A)\cos(t
\ome_2 B)\cos(t\ome_3 C)d\ome\bigg]\ .\cr}\eqno (1.2)$$
As a simple illustration, we derive the well-known formulas for the
solution of the initial value problem for the wave equation in $R^n$
from the elementary solution for the one-dimensional case (Example~4.1).
Note that we do not use explicitly the spherical symmetry of $\Del$. 
We also solve the initial value problem for the Klein-Gordon equation.
\par
In Section~3 we derive formulas for $\csABt$ when $A$ and $B$ do not
commute. That the expression should involve a limit follows from the
Trotter product formula for $e^{-(A^2+B^2)}$. Formal ``inversion'' of 
the Trotter formula leads to a limit formula for $\csABt$. This formula 
may be justified weakly (Proposition~3.1) or strongly on appropriate 
analytic vectors \cite{N}\ (Theorem~3.2). One could express $\cos\left(
\sqrt{A_1^2+\cdots+A_n^2}\,t\right)$ for general $n$ (Remark~3.4) but 
the resulting formulas are quite heavy. The formulas (such as (3.4)) are 
reminiscent of those appearing in path integrals (such as the 
Feynman-Kac formulas for parabolic or Schr\"odinger equations). This 
similarity consists of the integrations being extended over balls of 
increasing dimensions; due to the hyperbolic nature of the problem, we
have here in addition also differentiations to high orders.\par
To illustrate the formula for the non-commutative case we compute the 
propagators of the harmonic oscillator
$${\partial^2 u\over\partial t^2}={\partial^2 u\over\partial x^2}-x^2u$$
(Example~4.2) and of simple hypoelliptic sum of squares operators,
such as
$${\partial^2 u\over\partial t^2}={\partial^2 u\over\partial x_1^2}+
x_1^2{\partial^2 u\over\partial x_2^2}$$
and
$${\partial^2 u\over\partial t^2}=\Delta_H u$$
where $\Delta_H$ is the Laplacian on the Heisenberg group (Example~4.3).
The expressions obtained are not very explicit.\par
Note that the operator ${\sin\left(\sqrt{A^2+B^2}\,t\right)\over
\sqrt{A^2+B^2}}$, and more generally, the operator\linebreak 
$\sin\left(\sqrt{A_1^2+\cdots+A_n^2}\,t\right)\over\sqrt{A_1^2+\cdots+
A_n^2}$ may be represented by similar formulas, the only difference 
being dropping the left-most $\ppt$ from the corresponding formulas 
for\linebreak 
$\cos\left(\sqrt{A^2+B^2}\,t\right)$ $\left(\cos\left(t\sqrt{A_1^2+
\cdots+A_n^2}\,\right)\right)$. The proofs are similar to those of (2.2) 
or of Theorem~3.1.\par
I am very much indebted to V.~Katsnelson for simulating discussions.
\par
\section{The commutative case}
\noindent {\bf Theorem~2.1}\quad {\em Let\/ $H$ be a Hilbert space, and 
let the\/} ({\em possibly unbounded\/}) {\em self-adjoint operators\/ 
$A_1,\ldots,A_n$ commute\/}.\hfill\break
(i) {\em If\/ $n$ is odd, $n=2m+1$, then}
$$\eqalign{&\cos\left(t\sqrt{A_1^2+\cdots+A_n^2}\,\right)={1\over 2(2
\pi)^m}\ppt\ltptr^{m-1}\bigg[t^{2m-1}\cr
\noalign{\vskip 5pt}
&\qquad\times\int_{S^{2m}}\cos(t\ome_1 A_1)\cdots(t\ome_n 
A_n)d\ome\bigg]\cr}\eqno (2.1)$$
({\em here\/ $d\ome$ denotes the surface measure on\/} $S^{2m}$).\hfill
\break
(ii) {\em If\/ $n$ is even, $n=2m$, then}
$$\eqalign{&\cos\left(t\sqrt{A_1^2+\cdots+A_n^2}\,\right)={1\over (2
\pi)}^m\ppt\ltptr^{m-1}\bigg[t^{2m-1}\cr
\noalign{\vskip 5pt}
&\qquad\times\int_{\ome_1^2+\cdots+\ome_{2m}^2\le 1} {\cos(t\ome_1 A_1)
\cdots\cos(t\ome_{2m} A_{2m})\over\sqrt{1-(\ome_1^2+\cdots+
\ome_{2m}^2)}}\ d\ome_1\cdots d\ome_{2m}\bigg]\ .\cr}\eqno (2.2)$$\par
\noindent
{\bf Proof.}\quad Set $L=\sqrt{\sum_{i=1}^n A_i^2}\,$. Note that $L$
is a non-negative self-adjoint operator in $H$. To prove (i), note
that by commutativity we have for $\rho>0$ that 
$$e^{-L^2\rho}=\prod_{i=1}^n e^{-A_i^2\rho}\ .\eqno (2.3)$$
The well-known transmutation formula (e.g.\ \cite{K})
$$e^{-B^2\rho}={1\over\sqrt{4\pi\rho}}\iii e^{-{t^2\over 4\rho}}\cos(Bt)
dt\ ,\eqno (2.4)$$
valid for $B$ self-adjoint, implies that 
$$e^{-B^2\rho}={1\over\sqrt{\pi\rho}}\izi e^{-{t^2\over 4\rho}}\cos(Bt)
dt\ ,\eqno (2.4^{\prime})$$
Applying (2.4) to each $A_i$, we find that
$$\eqalign{\prod_{i=1}^n e^{-A_i^2\rho}&={1\over(4\pi\rho)^{n/2}}\iii
e^{-{t_1^2\over 4\rho}}\cos(A_1t_1)dt_1\iii e^{-{t_2^2\over 4\rho}}
\cos(A_2t_2)dt_2\cdots\cr
\noalign{\vskip 5pt}
&\qquad\times\iii e^{-{t_n^2\over 4\rho}}\cos(A_nt_n)dt_n\cr
\noalign{\vskip 5pt}
&={1\over(4\pi\rho)^{n/2}}\int_{R^n}e^{-{t_1^2+\cdots+t_n^2\over 4\rho}}
\cos(A_1t_1)\cdots\cos(A_nt_n)dt\cr
\noalign{\vskip 5pt}
&={1\over(4\pi\rho)^{n/2}}\izi t^{n-1}e^{-{t^2\over 4\rho}}\int_{S^{n-1}}
\cos(t\ome_1 A_1)\cdots\cos(t\ome_n A_n)d\ome dt\ .\cr}\eqno (2.5)$$
Setting $B=L$ in (2.4$^{\prime}$) we see from (2.3) and (2.5) that
$$\eqalign{{1\over\sqrt{\pi\rho}}\izi e^{-{t^2\over 4\rho}}\cos(Lt)dt=&
{1\over(4\pi\rho)^{n/2}}\int_0^{\infty} t^{n-1}e^{-{t^2\over 4\rho}}\cr
\noalign{\vskip 5pt}
&\times\int_{S^{n-1}}\cos(t\ome_1 A_1)\cdots\cos(t\ome_n A_n)d\ome 
dt\cr}$$
so that
$$\izi e^{-{t^2\over 4\rho}}\cos(Lt)dt={1\over 2(4\pi\rho)^m}\izi 
e^{-{t^2\over 4\rho}} t^{2m}F(t;A_1,\ldots,A_n)dt\eqno (2.6)$$
where
$$F(t;A_1,\ldots,A_n)=\int_{S^{2m}}\cos(t\ome_1 A_1)\cos(t\ome_2A_2)
\cdots\cos(t\ome_n A_n)d\ome\ .$$
We apply successive integrations by parts to eliminate the factor
$\rho^{-m}$ that appears in the right hand side of (2.6). Thus,
$$\eqalign{&{1\over(4\pi\rho)^m}\izi e^{-{t^2\over 4\rho}} t^{2m} F(t;
A_1,\ldots,A_n)dt\cr
\noalign{\vskip 5pt}
&\qquad =-{1\over 2\pi}{1\over(4\pi\rho)^{m-1}}\izi\ppt\left(e^{-{t^2
\over 4\rho}}\right)t^{2m-1}F(t;A_1,\ldots,A_n)dt\cr
\noalign{\vskip 5pt}
&\qquad ={1\over 2\pi}{1\over(4\pi\rho)^{m-1}}\izi e^{-{t^2\over 4\rho}}
\ppt\left[t^{2m-1}F(t;A_1,\ldots,A_n)\right]dt\cr
\noalign{\vskip 5pt}
&\qquad =-{1\over(2\pi)^2}{1\over(4\pi\rho)^{m-2}}\izi\ltptr\left(
e^{-{t^2\over 4\rho}}\right)\ppt\left[t^{2m-1}F(t;A_1,\ldots,A_n)\right]
dt\cr
\noalign{\vskip 5pt}
&\qquad ={1\over(2\pi)^2}{1\over(4\pi\rho)^{m-2}}\izi e^{-{t^2\over 4
\rho}}\left[\ppt\ltptr\right]\left(t^{2m-1}F(t;A_1,
\ldots,A_n)\right]dt\cr
\noalign{\vskip 5pt}
&\qquad =\cdots={1\over(2\pi)^m}\izi e^{-{t^2\over 4\rho}}\ppt
\ltptr^{m-1}\left[t^{2m-1}F(t;A_1,\ldots,A_n)\right]dt\cr}$$
Substituting in (2.6) and applying the uniqueness theorem for the
Laplace transform, we obtain (2.1).\par
One could derive (2.2) from (2.1) by setting $A_{2m+1}=0$ (essentially
a method of descent). We prefer an alternative method. Let $a_1,\ldots, 
a_{2m}$ be scalars. Then (using multi-index notation) 
$$\eqalign{&\int_{\pll\ome\pll\le 1}{\cos(\ome_1 a_1t)\cdots\cos(
\ome_{2m}a_{2m}t)\over\sqrt{1-\pll\ome\pll^2}} d\ome_1\cdots d
\ome_{2m}\cr
\noalign{\vskip 5pt}
&\qquad =\sum_{\alp}(-1)^{\alp}{t^{2\vert\alp\vert} a^{2\alp}\over 
(2\alp)!}\int_{\pll\ome\pll\le 1}{\ome^{2\alp}\over\sqrt{1-\pll\ome
\pll^2}} d\ome_1\cdots d\ome_{2m}\ .\cr}\eqno (2.7)$$
The integral on the right hand side of (2.7) is essentially the
Dirichlet integral \cite{W}. In fact, 
$$\int_{\pll\ome\pll^2\le 1}{\ome^{2\alp}\over\sqrt{1-\pll\ome\pll^2}}
d\ome_1\cdots d\ome_{2m}=2^{2m}\int_{\ome\ge 0,\,\pll\ome\pll^2\le 1}
{w^{2\alp}\over\sqrt{1-\pll\ome\pll^2}}d\ome_1\cdots d\ome_n\ .$$
Set $z_i=\ome_i^2$, $1\le i\le 2m$. Then the integral on the right hand 
side of (2.7) is equal to
$$\int_{z\ge 0,\,z_1+\cdots z_{2m}\le 1}z_1^{\alp_1-{1\over 2}}\cdots 
z_{2m}^{\alp_{2m}-{1\over 2}}(1-z_1-\ldots-z_{2m})^{-1/2}dz_1\cdots 
dz_{2m}\ .$$
According to formulas (7.7.4) and (7.7.5) in \cite{W}, the value of
this integral is 
$${\Gamma\left(\alp_1+{1\over 2}\right)\cdots\Gamma\left(\alp_{2m}+
{1\over 2}\right)\Gamma\left({1\over 2}\right)\over\Gamma\left(\vert
\alp\vert+m+{1\over 2}\right)}\ .$$
Substituting in (2.7) and recalling the duplication formula for the
Gamma function
$$\Gamma\left(k+{1\over 2}\right)={\sqrt{\pi}\,\Gamma(2k)\over 2^{2k-1} 
\Gamma(k)}$$
we see that
$$\eqalign{&\int_{\pll\ome\pll^2\le 1}{\cos(\ome_1 a_1t)\cdots\cos(
\ome_{2m}a_{2m}t)\over\sqrt{1-\pll\ome\pll^2}} d\ome_1\cdots d\ome_{2m}\cr
\noalign{\vskip 5pt}
&=\sum_{\alp}(-1)^{\vert\alp\vert}{t^{2\vert\alp\vert}a^{2\alp}\over(2
\alp)!}{\Gamma\left(\alp_1+{1\over 2}\right)\cdots\Gamma\left(\alp_{2m}
+{1\over 2}\right)\Gamma\left({1\over 2}\right)\over\Gamma\left(\vert
\alp\vert+m+{1\over 2}\right)}\cr
\noalign{\vskip 5pt}
&=\sum_{\alp}{(-1)^{\vert\alp\vert}\over(2\alp)!}{t^{2\vert\alp\vert}
a^{2\alp}\pi^m\over 2^{2\vert\alp\vert-2m}}{\Gamma(2\alp_1)\cdots\Gamma
(2\alp_{2m})2^{2\vert\alp\vert+2m-1}\Gamma(\vert\alp\vert+m)\over\Gamma(
\alp_1)\cdots\Gamma(\alp_{2m})\Gamma(2\vert\alp\vert+2m)}\cr
\noalign{\vskip 5pt}
&=\sum_{\alp}{(-1)^{\vert\alp\vert}t^{2\vert\alp\vert}a^{2\alp}\pi^m 
2^{4m}\over(2\alp_1)\cdots(2\alp_{2m})\Gamma(\alp_1)\cdots\Gamma(
\alp_{2m})}{1\over 2^{\vert\alp\vert+m}(2m+2\vert\alp\vert-1)}\cr
\noalign{\vskip 5pt}
&\qquad\times{1\over(2m+2\vert\alp\vert-3)\cdots 1}\cr
\noalign{\vskip 5pt} 
&=\sum_{\alp}{(-1)^{\vert\alp\vert}t^{2\vert\alp\vert}a^{2\alp}\over
\alp!}{\pi^m 2^{4m}\over 2^{2m}}{1\over 2^{\vert\alp\vert+m}}{1\over(2m
+2\vert\alp\vert-1)(2m+2\vert\alp\vert-3)\cdots 1}\cr
\noalign{\vskip 5pt}
&=\sum_{\alp}{(-1)^{\vert\alp\vert}t^{2\vert\alp\vert}a^{2\alp}(2\pi)^m
\over2^{\vert\alp\vert}\alp!}{1\over(2m+2\vert\alp\vert-1)(2m+2\vert\alp
\vert-3)\cdots 1}\ .\cr}$$
Hence
$$\eqalign{&t^{2m-1}\int_{\pll\ome\pll^2\le 1}{\cos(\ome_1 a_1t)\cdots
\cos(\ome_{2m} a_{2m}t)\over\sqrt{1-\pll\ome\pll^2}}d\ome_1\cdots d
\ome_{2m}\cr
\noalign{\vskip 5pt}
&\qquad =\sum_{\alp}{(-1)^{\vert\alp\vert}(2\pi)^m\over2^{\vert\alp
\vert}\alp!}{t^{2\vert\alp\vert+2m-1}a^{2\alp}\over(2\vert\alp\vert+2m-
1)(2\vert\alp\vert+2m-3)\cdots 1}\ .\cr}\eqno (2.8)$$
But 
$$\eqalign{\ltptr t^{2\vert\alp\vert+2m-1}&=(2\vert\alp\vert+2m-1)t^{2
\vert\alp\vert+2m-3},\cdots,\ltptr^{m-1}t^{2\vert\alp\vert+2m-1}\cr
\noalign{\vskip 5pt}
&=(2\vert\alp\vert+2m-1)\cdots(2\vert\alp\vert+3)t^{2\vert\alp\vert+1}
\ .\cr}$$
Applying to (2.8), we get
$$\eqalign{&\ltptr^{m-1}\bigg[t^{2m-1}\int_{\pll\ome\pll^2\le 1}{\cos(
\ome_1 a_1t)\cdots\cos(\ome_{2m} a_{2m}t)\over\sqrt{1-\pll\ome\pll^2}}
d\ome_1\cdots d\ome_{2m}\bigg]\cr
\noalign{\vskip 5pt}
&\qquad =\sum_{\alp}{(-1)^{\vert\alp\vert}(2\pi)^m a^{2\alp}\over 
2^{\vert\alp\vert}\alp!}{t^{2\vert\alp\vert+1}\over(2\vert\alp\vert+1)
(2\vert\alp\vert-1)\cdots 1}\cr}$$
and
$$\eqalign{&{1\over(2\pi)^m}\ppt\ltptr^{m-1}\bigg[t^{2m-1}\int_{\pll\ome
\pll^2\le 1}{\cos(\ome_1 a_1t)\cdots\cos(\ome_{2m} a_{2m}t)\over\sqrt{1-
\pll\ome\pll^2}}d\ome_1\cdots d\ome_{2m}\bigg]\cr
\noalign{\vskip 5pt}
&\qquad =\sum_{\alp}{(-1)^{\vert\alp\vert}a^{2\alp}t^{2\vert\alp\vert}
\over 2^{\vert\alp\vert}\alp!(2\vert\alp\vert-1)\cdots 1}\ .\cr}\eqno 
(2.9)$$
On the other hand, for every non-negative integer $k$,
$$\sum_{\vert\alp\vert = k}{a^{2\alp}\over\alp!}={1\over k!}\sum_{\vert\alp
\vert = k}{k!\over\alp!}a^{2\alp}={1\over k!}(a_1^2+\cdots+a_{2m}^2)^k
\ .$$
Inserting in (2.9) we conclude that
$$\eqalign{&{1\over(2\pi)^m}\ppt\ltptr^{m-1}\bigg[t^{2m-1}\int_{\pll\ome
\pll^2\le 1}{\cos(\ome_1 a_1t)\cdots\cos(\ome_{2m} a_{2m}t)\over\sqrt{1-
\pll\ome\pll^2}}d\ome_1\cdots d\ome_{2m}\bigg]\cr
\noalign{\vskip 5pt}
&\qquad =\sum_{k=0}^{\infty}{(-1)^k\over 2^kk!}{t^{2k}(a_1^2+\cdots+
a_{2m}^2)^k\over(2k-1)\cdots}=\sum_{k=1}^{\infty}{(-1)^k(\vert a\vert
t)^{2k}\over(2k)!}=\cos(\vert a\vert t)\cr}$$
proving (2.2) for the case $A_i=a_i$, $i=1,\ldots,2m$. In general, let
$dE^{(i)}(\lam)$ denote the spectral measure of $A_i$, $1\le i\le 2m$. 
The commutativity assumption implies that
$$\eqalign{&\cos\left(t\sqrt{A_1^2+\cdots+A_{2m}^2}\,\right)\cr
\noalign{\vskip 5pt}
&\qquad =\int\cos\left(t\sqrt{\lam_1^2+\cdots\lam_{2m}^2}\,\right)
dE^{(1)}(\lam_1)\cdots dE^{(2m)}(\lam_{2m})\ .\cr}\eqno (2.10)$$
Application of (2.2) (case of real numbers) to (2.10) implies that 
$$\eqalign{&\cos\left(t\sqrt{A_1^2+\cdots+A_{2m}^2}\,\right)=\int 
{1\over(2\pi)^m}\ppt\ltptr^{2m-1}\bigg[t^{2m-1}\cr
\noalign{\vskip 5pt}
&\qquad\times\int_{\pll\ome\pll^2\le 1}{\cos(\ome_1\lam_1 t)\cdots\cos(
\ome_{2m}\lam_{2m} t)\over\sqrt{1-\pll\ome\pll^2}} d\ome_1\cdots d
\ome_{2m}\bigg]\cr
\noalign{\vskip 5pt}
&\qquad\times dE^{(1)}(\lam_1)\cdots dE^{(2m)}(\lam_{2m})\cr
\noalign{\vskip 5pt}
&={1\over(2\pi)^m}\ppt\ltptr^{2m-1}\bigg\{t^{2m-1}\cr
\noalign{\vskip 5pt}
&\qquad\times\int_{\pll\ome\pll^2\le 1}\bigg[\int{\cos(\ome_1\lam_1 t)
\cdots\cos(\ome_{2m}\lam_{2m} t)\over\sqrt{1-\pll\ome\pll^2}}\cr
\noalign{\vskip 5pt}
&\qquad\times dE^{(1)}(\lam_1)\cdots dE^{(2m)}(\lam_{2m})\bigg] d\ome_1
\cdots d\ome_{2m}\bigg\}\cr
\noalign{\vskip 5pt}
&={1\over(2\pi)^m}\ppt\ltptr^{2m-1}\bigg[t^{2m-1}\int{\cos(\ome_1
A_1t) \cdots\cos(\ome_{2m}
A_{2m}t)\over\sqrt{1-\pll\ome\pll^2}}d\ome_1\cdots d\ome_{2m}\bigg]\cr}$$
proving (2.2).\smboc\par
\section{The non-commutative case}
Let $A$, $B$ be (not necessarily commuting) self-adjoint operators on
the Hilbert space $H$. Then the operators $A^2$, $B^2$ are semi-bounded.
Recall the Trotter product formula \cite{RS}
$$e^{-\rho(A^2+B^2)}=\lim_{m\rar\infty}\left[e^{-\rho{A^2\over m}}
e^{-\rho{B^2\over m}}\right]^m\eqno (3.1)$$
valid for $\rho>0$ (or, more generally, for $\Ree\rho\ge 0$). Setting
$n=2m+1$,
$$\eqalign{&A_1={A\over\sqrt{m}}\ ,\ A_2={B\over\sqrt{m}}\ ,\ A_3={A\over
\sqrt{m}}\ ,\ \ldots\ ,\ A_{2m-1}={A\over\sqrt{m}}\ ,\ A_{2m}={B\over
\sqrt{m}}\ ,\cr 
\noalign{\vskip 5pt}
&A_{2m+1}=0\ ,\cr}$$
we infer from (2.5) that
$$\eqalign{&\left[e^{-\rho{A^2\over m}}e^{-\rho{B^2\over m}}\right]^m=
{1\over(4\pi\rho)^{m+1/2}}\izi t^{n-1}e^{-{t^2\over 4\rho}}\int_{S^{2m}}
\cos\left({t\ome_1 A\over\sqrt{m}}\right)\cr
\noalign{\vskip 5pt}
&\qquad\times\cos\left({t\ome_2 B\over\sqrt{m}}\right)\cdots\cos\left(
{t\ome_{2m-1} A\over\sqrt{m}}\right)\cos\left({t\ome_{2m} B\over
\sqrt{m}}\right)d\ome dt\ .\cr}\eqno (3.2)$$
Note that $d\ome={2d\ome_1\cdots d\ome_{2m}\over\sqrt{1-(\ome_1^2+\cdots
+\ome_{2m}^2)}}\,$. Integrating by parts as in Section~2, we see from
(3.1) and (3.2) that
$$\eqalign{&\izi e^{-{t^2\over 4\rho}}\csABt dt=\lim_{m\rar\infty}
\int_0 e^{-{t^2\over 4\rho}}{1\over(2\pi)^m}\ppt\ltptr^{m-1}
\bigg[t^{2m-1}\cr
\noalign{\vskip 5pt}
&\qquad\times\int_{\ome_1^2+\cdots+\ome_{2m}^2\le 1}\cr
\noalign{\vskip 5pt}
&\qquad\times {\cos\left({t\ome_1 A\over\sqrt{m}}\right)\cos\left(
{t\ome_2 B\over\sqrt{m}}\right)\cdots\cos\left({t\ome_{2m-1} A\over
\sqrt{m}}\right)\cos\left({t\ome_{2m} B\over\sqrt{m}}\right)\over\sqrt{1-
(\ome_1^2+\cdots+\ome_{2m}^2)}} d\ome_1\cdots d\ome_{2m}\bigg]dt\ .\cr}
\eqno (3.3)$$
The relation (3.3) suggests, at least formally, that
$$\eqalign{&\csABt=\lim_{m\rar\infty}{1\over(2\pi)^m}\ppt\ltptr^{m-1}
\bigg[t^{2m-1}\int_{\ome_1^2+\cdots+\ome_{2m}^2\le 1}\cr
\noalign{\vskip 5pt}
&\qquad\times {\cos\left({t\ome_1 A\over\sqrt{m}}\right)\cos\left(
{t\ome_2 B\over\sqrt{m}}\right)\cdots\cos\left({t\ome_{2m-1} A\over
\sqrt{m}}\right)\cos\left({t\ome_{2m} B\over\sqrt{m}}\right)\over\sqrt{1-
(\ome_1^2+\cdots+\ome_{2m}^2)}} d\ome_1\cdots d\ome_{2m}\bigg]\ .\cr}
\eqno (3.4)$$
In the rest of this section we try to interpret and to justify the
heuristically obtained formula (3.4). Note that for each $m$, the
function 
$$C_m(t)=\int_{\ome_1^2+\cdots+\ome_{2m}^2\le 1}\cos\left({t\ome_1 A
\over\sqrt{m}}\right)\cdots\cos\left({t\ome_{2m} B\over\sqrt{m}}\right)
d\ome_1\cdots d\ome_{2m}$$
is a holomorphic function from the complex numbers to the space of
bounded operators on $H$. Thus we ``only'' have to make precise the
convergence in (3.4). Our first result is about weak convergence.\par
\medskip\noindent
{\bf Proposition~3.1}\quad {\em Let\/ $A^2+B^2$ be essentially 
self-adjoint on\/ $D(A^2)\cap D(B^2)$. For every positive integer\/ $m$ 
and\/ $h\in H$, define the function\/ $F_m(t;A,B;h)\in C(R^1;H)$ by}
$$\eqalign{&F_m(t;A,B;h)={1\over(2\pi)^m}\ppt\ltptr^{m-1}
\bigg[t^{2m-1}\cr
\noalign{\vskip 5pt}
&\qquad\times\int_{\ome_1^2+\cdots+\ome_{2m}^2\le 1}{\cos\left({\ome_1
tA\over\sqrt{m}}\right)\cdots\cos\left({\ome_{2m} tB\over\sqrt{m}}
\right)\over\sqrt{1-(\ome_1^2+\cdots+\ome_{2m}^2)}} d\ome_1\cdots 
d\ome_{2m}h\bigg]\ .\cr}\eqno (3.5)$$
{\em Then the functions\/ $F_m(t;A,B;h)$ converge as vector valued temperate 
distributions\/} ({\em i.e., in\/} $S^{\prime}(R^1;H)$) {\em to\/} 
$\cos\left(\sqrt{A^2+B^2}\,t\right)h$.\par
\medskip\noindent
{\bf Proof.}\quad Let $\Ree\rho\ge 0$, set ${1\over 2\rho}=i\sig$
so that $\Imm\sig\le 0$, and let ${t^2\over 2}=s$. It follows
from  (3.2) and (3.3) that
$$\left[e^{-{A^2\over 2i\sig m}}e^{-{B^2\over 2i\sig m}}\right]^mh=
\sqrt{2i\sig\over\pi}\izi F_m\left(\sqrt{2s}\,;A,B;h\right)e^{-is\sig}
{ds\over\sqrt{2s}}\eqno (3.6)$$
and, similarly,
$$e^{-{A^2+B^2\over 2i\sig}}h=\sqrt{2i\sig\over\pi}\izi\cos\left(
\sqrt{A^2+B^2}\,\sqrt{2s}\,\right)he^{-is\sig}{ds\over\sqrt{2s}}\eqno
(3.7)$$
when $\sqrt{2i\sig}$ is chosen to be positive for $\sig$ on the
negative imaginary axis. The left hand sides of (3.6) and (3.7) are
uniformly bounded (by $\pll h\pll$) for $\Imm\sig\ge 0$, and are
holomorphic for $\Imm\sig>0$. Hence $\sqrt{\pi\over 2i\sig}\left[
e^{-{A^2\over 2i\sig m}}e^{-{B^2\over 2i\sig m}}\right]^mh$ is the
Fourier transform of $F_m\left(\sqrt{2s}\,;A,B;h\right)/\sqrt{2s}$ and
$\sqrt{\pi\over 2i\sig}e^{-{A^2+B^2\over 2i\sig}}h$ is the Fourier
transform of ${\cos\left(\sqrt{A^2+B^2}\,\sqrt{2s}\,\right)\over
\sqrt{2s}}$. By the Trotter product formula, the left hand sides of
(3.6) converge to the left side of (3.7) for every $\sig\in R^1\backslash
\{0\}$. Let $\psi(\sig)\in S(R^2,H)$ be arbitrary. By the dominated
convergence theorem
$$\int\Bigl\langle\left(e^{-{A^2\over 2i\sig m}}e^{-{B^2\over 2i\sig m}}
\right)^mh,\psi(\sig)\Bigr\rangle {d\pi\over\sqrt{\sig}}\rar\int
\Bigl\langle e^{-{A^2+B^2\over 2i\sig}}h,\psi(\sig)\Bigr\rangle
{d\sig\over\sqrt{\sig}}\ .$$
The continuity of the inverse Fourier transform in $S^{\prime}(R^1,H)$
implies that the sequence 
$F_m\left(\sqrt{2s}\,;A,B;h\right)/\sqrt{2s}$ converges, in the 
temperate distribution sense, to 
$\cos\left(\sqrt{A^2+B^2}\,\sqrt{2s}\,\right)\sqrt{2s}$. Passing back 
to $t=\sqrt{2s}$ and noting that \linebreak $F_m(t;A,B;h)$ and $\cos\left(\sqrt{A^2
+B^2}\,t\right)$ are even functions, we obtain the result.\smboc\par
One could possibly get a more precise version of Proposition~3.1 by
using a strong version of the Trotter product formula and by estimating 
the convergence in an appropriate Sobolev norm. We state now a result
on pointwise convergence, under stronger assumptions on $h$.\par
\medskip\noindent
{\bf Theorem~3.1}\quad {\em Let\/ $h\in H$ be such that there exist 
constants\/ $C$, $K$ so that for every multi-index\/ $\alp$ with\/ $2m$ 
components\/ $h\in D\left(A^{\alp_1}B^{\alp_2}\cdots A^{\alp_{2n-1}}
B^{\alp_{2n}})\right)$ and}
$$\pll A^{\alp_1}B^{\alp_2}\cdots A^{\alp_{2n-1}}B^{\alp_{2n}}h\pll\le
CK^{\vert\alp\vert}\vert\alp\vert!\ .\eqno (3.8)$$
{\em Then\/ $F_m(t;A,B;h)\rar\cos\left(\sqrt{A^2+B^2}\,t\right)h$ for\/
$\vert t\vert<{1\over\sqrt{2}K}$, uniformly in compact subintervals of\/}
$\left(-{1\over\sqrt{2}K},{1\over\sqrt{2}K}\right)$.\par
\medskip\noindent
{\bf Proof.}\quad Note that for every positive integer $n$,
$$(A^2+B^2)^n=\sum_{\epsi_i=0,1,\ 1\le i\le n} A^{2\epsi_1}
B^{2(1-\epsi_1)}\cdots A^{2\epsi_n}B^{2(1-\epsi_n)}\ .\eqno (3.9)$$
By (3.8) we have for $h$ satisfying the assumptions that 
$$\pll(A^2+B^2)^nh\pll\le 2^nCK^{2n}(2n)!\eqno (3.10)$$
Expanding $\cos\left(\sqrt{A^2+B^2}\,t\right)h$ in Taylor series, we
see that for every $N$,
$$\eqalign{&\pll\cos\left(\sqrt{A^2+B^2}\,t\right)h-\sum_{n\le 
N}{(-1)^n\over (2n)!} t^{2n}(A^2+B^2)^nh\pll\cr
\noalign{\vskip 5pt}
&\qquad\le\sum_{n=N+1} {\vert t\vert^{2n}\over(2n)!} 2^nCK^{2n}(2n)!=
C\sum_{n=N+1}^{\infty}\left(\sqrt{2}\,\vert t\vert K\right)^{2n}=
{C\left(\sqrt{2}\,tK\right)^{2N+2}\over 1-2t^2K^2}\cr}\eqno (3.11)$$
or that the series converges absolutely and uniformly in compact
subsets of $\vert t\vert<{1\over\sqrt{2}\,K}$.\par
Expanding similarly $\cos\left({\ome_1 tA\over\sqrt{m}}\right)\cos
\left({\ome_2 tB\over\sqrt{m}}\right)\cdots\cos\left({\ome_{2m} tB\over
\sqrt{m}}\right)h$ in Taylor series, we get the estimate 
$$\eqalign{&\Bigl\vert\cos\left({\ome_1 tA\over\sqrt{m}}\right)\cos
\left({\ome_2 tB\over\sqrt{m}}\right)\cdots\cos\left({\ome_{2m-1} tA
\over\sqrt{m}}\right)\cos\left({\ome_{2m} tB\over\sqrt{m}}\right)h-
\sum_{n=0}^N(-1)^n\cr
\noalign{\vskip 5pt}
&\qquad\times {t^{2n}\over m^n}\sum_{\alp\in Z_+^{2m}\atop\vert\alp\vert
=n}{(\ome_1 A)^{2\alp_1}(\ome_2 B)^{2\alp_2}\cdots(\ome_{2n-1}A)^{2
\alp_{2m-1}}(\ome_{2m}B)^{2\alp_{2m}}h\over(2\alp)!}\Bigr\vert\cr
\noalign{\vskip 5pt}
&\quad\le\sum_{n=N+1}^{\infty} {\vert t\vert^{2n}\over m^n}\sum_{\alp
\in Z_+^{2m}\atop\vert\alp\vert=n}{\pll A^{2\alp_1}B^{2\alp_2}\cdots 
A^{2\alp_{2m-1}}B^{2\alp_{2m}}h\pll\over(2\alp)!}\cr
\noalign{\vskip 5pt}
&\quad\le\sum_{n=N+1}^{\infty} {t^{2n}\over m^n}\sum_{\alp\in Z_+^{2m}
\atop\vert\alp\vert=n}{CK^{2n}(2n)!\over(2\alp)!}\le C\sum_{n=N+1}^{\infty}
(tK)^{2n}{(2m)^{2n}\over m^n}\cr
\noalign{\vskip 5pt}
&\quad =C\sum_{n=N+1}^{\infty} (4t^2K^2m)^n\cr}$$
so that this expansion converges absolutely and uniformly in the
product of $\pll\ome\pll\le 1$ with compact subsets of $\vert t\vert
<{1\over 2K\sqrt{m}}$; ${1\over 2K\sqrt{m}}$ is positive for every 
{\em fixed\/} $m$. Hence
we may perform the integration in (3.5) term by term. Evaluating the
integral $\int_{\pll\ome\pll\le 1}{\ome^{2\alp}\over\sqrt{1-\pll\ome
\pll^2}} d\ome_1\cdots d\ome_{2m}$ as in Section~2, we obtain the
non-commutative version of (2.9) with $a_i={A\over\sqrt{m}}$ for $i$
even, $a_i={B\over\sqrt{m}}$ for $i$ odd, i.e., that
$$\eqalign{F_m(t;A,B;h)=~&\sum_{n=0}^{\infty} {(-1)^n\over 2^n}{t^{2n}
\over m^n(2n-1)\cdots 1}\cr
\noalign{\vskip 5pt}
&\times\sum_{\alp\in Z_+^{2m}\atop\vert\alp\vert=n}{A^{2\alp_1}
B^{2\alp_2}\cdots A^{2\alp_{2m-1}}B^{2\alp_{2m}}h\over\alp!}\ , \cr}\eqno 
(3.12)$$
the series converging absolutely and uniformly in compact subsets of
$\vert t\vert<{1\over 2K\sqrt{m}}$. Noting that $2^nn!(2n-1)\cdots 1=(2n)!$,
we may write (3.12) in the form
$$F_m(t;A,B;h)=\sum_{n=0}^{\infty} {(-1)^nt^{2n}\over(2n)!m^n}
\sum_{\alp\in Z_+^{2m}\atop\vert\alp\vert=n} {n!\over\alp!}A^{2\alp_1}
B_{2\alp_2}\cdots A^{2\alp_{2m-1}}B^{2\alp_m}h\ .\eqno (3.13)$$
The $n^{th}$ term (sum) in (3.11) is majorized according to the
assumption (3.8) by
$$\eqalign{&{\vert t\vert^{2n}\over(2n)!m^n}\sum_{\alp\in Z_+^{2m}\atop
\vert\alp\vert=n} {n!\over\alp!} CK^{2n}(2n)!\cr
\noalign{\vskip 5pt}
&\qquad =C{(\vert t\vert^2K^2)^n\over m^n}\sum_{\alp\in Z_+^{2m}\atop
\vert\alp\vert=n} {n!\over\alp!}={C(\vert t\vert^2K^2)^n\over m^n}(2m)^n
=C(2\vert t\vert^2K^2)^n\cr}$$
so that the radius of convergence of the series in (3.11) is ${1\over
\sqrt{2}\,K}$, and by analyticity of $F_m(t;A,B;h)$ (as a function of
$t$) the series expansion (3.13) is valid for $\vert t\vert<{1\over
\sqrt{2}\,K}$, {\em uniformly in\/} $m$.\par
Let $X$, $Y$ be non-commuting finite dimensional matrices. In this
case the convergence in the Trotter product formula \cite{RS}
$$e^{z(X+Y)}=\lim_{n\rar\infty}\left(e^{zX/m}e^{zY/m}\right)^m$$
is uniform in compact subset of the complex plane. Hence the Taylor
coefficients of the right hand side converge to those of the left hand
side, or
$${(X+Y)^n\over n!}=\lim_{m\rar\infty}\sum_{\alp\in Z_+^{2m}\atop\vert
\alp\vert=n}{X^{\alp_1}Y^{\alp_2}\cdots X^{\alp_{2n-1}}Y^{\alp_{2n}}
\over m^n\alp!}\eqno (3.14)$$
for every (fixed) $n$. It follows from (3.8) that for every (fixed)
$n$ and $h$,
$${(A^2+B^2)^n\over n!}h=\lim_{m\rar\infty}\sum_{\alp\in Z_+^{2m}\atop
\vert\alp\vert=n}{A^{2\alp_1}B^{2\alp_2}\cdots A^{2\alp_{2m-1}}
B^{\alp_{m}}\over m^n\alp!}h\ .\eqno (3.15)$$
For any given $\epsi=0$ and $T<{1\over\sqrt{2}\,K}$, choose $N$ such
that the right hand side of (3.11) and
$$\pll F_m(t;A,B;h)-\sum_{n=0}^N {(-1)^n\over(2n)!} {t^{2n}\over m^n}
\sum_{\alp\in Z_+^{2m}\atop\vert\alp\vert=n} {n!\over\alp!}A^{2\alp_1}
B_{2\alp_2}\cdots B^{2\alp_{2m}}h\pll$$
are each less than ${\epsi\over 3}$ for $\vert t\vert\le T$. By (3.15)
there exists a constant $M>0$ such that for every $0\le n\le N$ we have
$$\pll {t^{2n}\over(2n)!m^n}\sum_{\alp\in Z_+^{2m}\atop\vert\alp\vert=n}
{n!A^{2\alp_1}B^{2\alp_2}\cdots B^{2\alp_{2m}}\over\alp!}h-{t^{2n}\over
(2n)!}(A^2+B^2)^nh\pll<{\epsi\over 3(N+1)}$$
for $\vert t\vert\le T$ and $m>M$. Hence
$$\pll F_m(t;A,B;h)-\cos\left(\sqrt{A^2+B^2}\,t\right)h\pll<\epsi$$
for $\vert t\vert\le T$ and $m>M$.\smboc\par
\medskip\noindent
{\bf Remark~3.3}\quad If for every $K>0$ there exists a $C>0$ such that 
(3.8) holds for all $\alp$ then $F_m(t;A,B;h)\rar\cos\left(
\sqrt{A^2+B^2}\,t\right)h$ uniformly in compact subsets of the complex
$t$ plane. In particular, if $A$ and $B$ are bounded, then this holds
for every $h$ in $H$ (uniformly in bounded sets of $H$).\par
\medskip\noindent
{\bf Remark~3.4}\quad The Trotter product formula for any number of
non-commuting operators is well-known \cite{S}. In analogy to (3.4),
we may obtain the formula
$$\eqalign{&\cos\left(\sqrt{A_1^2+\cdots+A_q^2}\,t\right)=\lim_{m\rar 
\infty}{1\over(2\pi)^{mq/2}}\ppt\ltptr^{{mq\over 2}-1}\bigg[t^{mq-1}\cr
\noalign{\vskip 5pt}
&\qquad\times\int_{\ome_1^2+\cdots+\ome_{mq}^2\le 1}\cos
\left({t\ome_1 A_1\over\sqrt{m}}\right)\cdots\cos\left({t\ome_q A_q\over
\sqrt{m}}\right)\cos\left({t\ome_{q+1} A_1\over\sqrt{m}}\right)\cdots\cr
\noalign{\vskip 5pt}
&\qquad\times \cos\left({t\ome_{2q} A_q\over\sqrt{m}}\right){\cos\left(
{t\ome_{(m-1)q+1} A_1\over\sqrt{m}}\right)\cdots\cos\left({t\ome_{mq} 
A_q\over\sqrt{m}}\right)\over\sqrt{1-(\ome_1^2+\cdots+\ome_{mq}^2)}} d
\ome_1\cdots d\ome_{mq}\bigg]\cr}\eqno (3.16)$$
where $m$ run over all positive integers if $q$ is even, and over even
$m$ only if $q$ is odd. (Note that if $mq$ is odd then we may replace
the right hand side of (3.16) by a term involving integration over
$S^{mq}$, as in (2.1).) The formal limit (3.16) may be interpreted,
under suitable assumptions, as in Proposition~3.1 or in Theorem~3.2. We
leave the details to the diligent reader.\par
\medskip\noindent
{\bf Remark~3.5}\quad Analytic domination has been applied in the
context of Trotter products in \cite{A}.\par
\section{Illustrations and examples}
We observe that in all $\ome$-integrations appearing earlier we may
replace a term of the form $\cos(\ome_j T)$ by $\exp(i\ome_j T)$ (or
by $\exp(-i\ome_j T)$) since the rest of the integrand is even in
$\ome_j$ and $\sin(\ome_j T)$ is odd. We will use this observation
constantly in this section in order to simplify a few calculations.\par
\medskip\noindent
{\bf Example 4.1}\quad {\em The wave equation in\/ $R^n$ --- the method 
of ascent\/}.\par 
Recall that the operator ${1\over i}{d\over dx}$ has a (unique)
self-adjoint realization in $L^2(R^1)$ and the distribution kernels of
$\cos\left({t\over i}{d\over dx}\right)$, $\exp\left(it{1\over i}
{d\over dx}\right)$ and $\exp\left(-it{1\over i}{d\over dx}\right)$
(solutions of the one-dimensional wave equation) are given by
$$\eqalign{\cos\left({t\over i}{d\over dx}\right)(x,\ov{x})&={\del(x-
\ov{x}+t)+\del(x-\ov{x}-t)\over 2}\cr
\noalign{\vskip 5pt}
\exp\left(it{1\over i}{d\over dx}\right)(x,\ov{x})&=\del(x-\ov{x}+t)\cr
\noalign{\vskip 5pt}
\exp\left(-it{1\over i}{d\over dx}\right)(x,\ov{x})&=\del(x-\ov{x}-t)
\ .\cr}\eqno (4.1)$$
The simplest cases are $n=2,3$. In the two dimensional case we obtain
from (1.1) or (2.2) $(m=1)$ that
$$\eqalign{&\cos\left(t\sqrt{-{\ptl^2\over\ptl x^2}-{\ptl^2\over\ptl 
y^2}}\,\right)(x,y;\ov{x},\ov{y})\cr
\noalign{\vskip 5pt}
&\qquad ={1\over 2\pi}\ppt\bigg[t\int_{\ome_1^2+\ome_2^2\le 1}{\del(x-
\ov{x}+\ome_1 t)\del(y-\ov{y}+\ome_2 t)\over\sqrt{1-\ome_1^2-\ome_2^2}} 
d\ome_1 d\ome_2\bigg]\ .\cr}$$
But
$$\eqalign{&\int_{\ome_1^2+\ome_2^2\le 1}{\del(x-\ov{x}+\ome_1 t)\del(y-
\ov{y}+\ome_2 t)\over\sqrt{1-(\ome_1^2+\ome_2^2)}}d\ome_1 d\ome_2\cr
\noalign{\vskip 5pt}
&\qquad =\int_{\ome_1^2+\ome_2^2\le 1}{\del\left({x-\ov{x}\over t}+
\ome_1\right)\del\left({y-\ov{y}\over t}+\ome_2\right)\over t^2
\sqrt{1-(\ome_1^2+\ome_2^2)}}d\ome_1 d\ome_2\cr
\noalign{\vskip 5pt} 
&\qquad ={1\over t^2}{H\left(t^2-(x-\ov{x})^2-(y-\ov{y})^2\right)\over
\sqrt{1-\left[{(x-\ov{x})^2\over t^2}+{(y-\ov{y})^2\over t^2}\right]}}=
{1\over t}{H\left(t^2-(x-\ov{x})^2-(y-\ov{y})^2\right)\over\sqrt{t^2-(x
-\ov{x})^2+(y-\ov{y})^2}}\cr}$$
where $H$ is the Heaviside function. Hence
$$\cos\left(t\sqrt{-\left({\ptl\over\ptl x}\right)^2-\left({\ptl\over
\ptl y}\right)^2}\,\right)={1\over 2\pi}\ppt{H\left(t^2-(x-\ov{x})^2-
(y-\ov{y})^2\right)\over\sqrt{t^2-(x-\ov{x})^2+(y-\ov{y})^2}}\ .\eqno 
(4.2)$$
The three dimensional case is actually simpler:
$$\eqalign{&\cos\left(t\sqrt{-\left({\ptl\over\ptl x}\right)^2-{\ptl^2
\over\ptl y^2}-{\ptl^2\over\ptl z^2}}\,\right)(x,y,z;\ov{x},\ov{y},
\ov{z})\cr
\noalign{\vskip 5pt}
&\qquad ={1\over 4\pi}\ppt\bigg[t\int_{S^2}\left(e^{-it\ome_1 {1\over i}
{\ptl\over\ptl x}}\right)(x,\ov{x})\left(e^{-it\ome_2{1\over i}{\ptl
\over\ptl y}}\right)(y,\ov{y})\left(e^{-it\ome_3{1\over i}{\ptl\over
\ptl z}}\right)(z,\ov{z})d\ome\bigg]\cr
\noalign{\vskip 5pt}
&\qquad ={1\over 4\pi}\ppt\bigg[t\int_{S^2}\del(x-\ov{x}-t\ome_1)\del(y
-\ov{y}-t\ome_2)\del(z-\ov{z}-t\ome_3)d\ome\bigg]\cr
\noalign{\vskip 5pt}
&\qquad ={1\over 4\pi}\ppt\bigg[t\int_{S^2}\del\left(x-(\ov{x}+t\ome_1)
\right)\del\left(y-(\ov{y}+t\ome_2)\right)\del\left(z-(\ov{z}+t\ome_3)
\right)d\ome\bigg]\ .\cr}\eqno (4.3)$$
The formulas (4.2) and (4.3) are classical.\par
More generally, for $n\ge 3$ odd, $n=2m+1$, we have from (2.1) that
$$\eqalign{&\cos\left(t\sqrt{-\Del_n}\right)=\cos\left(
t\sqrt{-\sum_{j=1}^n {\ptl^2\over\ptl x_j^2}}\right)(x,\ov{x})\cr
\noalign{\vskip 5pt}
&\quad ={1\over 2(2\pi)^m}\ppt\ltptr^{m-1}\bigg[t^{2m-1}\cr
\noalign{\vskip 5pt}
&\qquad\times\int_{S^{2m}}\left(e^{-i\ome_1 t{1\over i}{\ptl\over\ptl 
x_1}}\right)(x_1,\ov{x}_1)\cdots\left(e^{-i\ome_n t{1\over i}{\ptl\over
\ptl x_n}}\right)(x_n,\ov{x}_n)d\ome\bigg]\cr
\noalign{\vskip 5pt}
&\quad ={1\over 2(2\pi)^m}\ppt\ltptr^{m-1}\bigg[t^{2m-1}\cr
&\qquad\times\int_{S^{2m}}\del (x_1-\ov{x}_1-\ome_1 t)\cdots\del(x_n-
\ov{x}_n-\ome_n t)d\ome\bigg]\ .\cr}\eqno (4.4)$$
But
$$1\cdot 3\cdots(n-2)\ome_n=1\cdot 3\cdots(n-2){2\pi^{n/2}\over\Gam
\left({n\over 2}\right)}=2(2\pi)^m\eqno (4.5)$$
and we recapture the classical formula \cite{F,H}. For general even
dimensions $n=2m$, we get from (2.2) that
$$\eqalign{&\cos\left(t\sqrt{-\Del_{2m}}\right)(x,\ov{x})={1\over
(2\pi)^m}\ppt\ltptr^{m-1}\bigg[t^{2m-1}\int_{\ome_1^2+\cdots+\ome_{2m}^2
\le 1}\cr
\noalign{\vskip 5pt}
&\qquad\times{\del(x_1-\ov{x}_1-\ome_1 t)\cdots\del(x_{2m}-\ov{x}_{2m}-
\ome_{2m} t)\over\sqrt{1-(\ome_1^2+\cdots+\ome_{2m}^2)}} d\ome_1\cdots
d\ome_{2m}\bigg]\ .\cr}\eqno (4.6)$$
Applying (4.5) (with $n=2m+1$) we see that the right hand side of (4.6) 
coincides with the well-known kernel of the operator mapping $u(x,0)$
to $u(x,t)$ (when ${\ptl u\over\ptl t}(x,0)=0$) for even $n$ \cite{F,H}.
Thus, we are able to build the solution of the initial value problem
for the $n$-dimensional wave equation from the (very elementary)
one-dimensional solution.\par
The Klein-Gordon operator $K$ is given by $K=\Del-a^2$, and we wish to
represent the solution of the initial value problem for the
Klein-Gordon equation
$${\ptl^2 u\over\ptl t^2}=\Del u-a^2u\eqno (4.7)$$
by calculating the distribution kernel of 
$$\cos\left(\sqrt{-K}\,t\right)=\cos\left(\sqrt{\sum_{j=1}^n\left(
{1\over i}{\ptl\over\ptl x_j}\right)^2+a^2}\,t\right)\ .$$ 
Set $A_j={1\over i}{\ptl\over\ptl x_j}$, $j=1,\ldots,n$, $A_{n+1}=a$. 
If $n$ is even, $n=2m$, we apply (2.1) to get
$$\eqalign{&\cos\left(\sqrt{-K}\,t\right)(x,\ov{x})={1\over 2(2\pi)^m}
\ppt\ltptr^{m-1}\bigg[t^{2m-1}\cr
\noalign{\vskip 5pt}
&\qquad\times\int_{S^{2m}}\cos(a\ome_{2m+1}t)\prod_{j=1}^{2m}\del(x_j-
\ov{x}_j-\ome_j t)d\ome\bigg]\cr
\noalign{\vskip 5pt}
&\quad ={1\over(2\pi)^m}\ppt\ltptr^{m-1}\bigg\{t^{2m-1}\int_{\ome_1^2+
\cdots+\ome_{2m-1}^2\le 1}{\cos\left(at\sqrt{1-\ome_1^2-\cdots-
\ome_{2m}^2}\,\right)\over\sqrt{1-\ome_1^2-\cdots-\ome_{2m}^2}}\cr
\noalign{\vskip 5pt}
&\qquad\times\prod_{j=1}^{2m}\left[\del\bigg({x_j-\ov{x}_j\over t}-
\ome_j\bigg)t^{-1}\right]d\ome_1\cdots d\ome_{2m}\bigg\}\cr
\noalign{\vskip 5pt}
&\quad ={1\over(2\pi)^m}\ppt\ltptr^{m-1}\bigg[{\cos\left(at\sqrt{1-
\sum_{j=1}^{2m}{(x_j-\ov{x}_j)^2\over t^2}}\,\right)H(t^2-\vert x-\ov{x}
\vert^2)\over t\sqrt{1-\sum_{j=1}^{2m}\left({x_j-\ov{x}_j\over t}
\right)^2}}\bigg]\ .\cr}$$
Hence
$$\eqalign{&\cos\left(\sqrt{-K}\,t\right)(x,\ov{x})\cr
\noalign{\vskip 5pt}
&\qquad ={1\over(2\pi)^m}\ppt\ltptr^{m-1}\bigg[{\cos\left(a\sqrt{t^2-
\vert x-\ov{x}\vert^2}\,\right)\over\sqrt{t^2-\vert x-\ov{x}\vert^2}}
H(t^2-\vert x-\ov{x}\vert^2)\bigg]\ .\cr}\eqno (4.8)$$
If $n$ is odd, $n=2m-1$, we get from (2.2) that
$$\eqalign{&\cos\left(\sqrt{-K}\,t\right)={1\over(2\pi)^m}\ppt
\ltptr^{m-1}\bigg[t^{2m-1}\cr
\noalign{\vskip 5pt}
&\qquad\times\int_{\ome_1^2+\cdots+\ome_{2m}^2\le 1}{\cos(a\ome_{2m}t)
\over\sqrt{1-(\ome_1^2+\cdots+\ome_{2m}^2)}}\cr
\noalign{\vskip 5pt}
&\qquad\times\prod_{j=1}^{2m-1}\del(x_j-\ov{x}_j-t\ome_i)d\ome_1\cdots
d\ome_{2m}\bigg]\cr
\noalign{\vskip 5pt}
&\quad ={1\over(2\pi)^m}\ppt\ltptr^{m-1}\bigg[\int_{_{-\sqrt{1-{\vert x-
\ov{x}\vert^2\over t^2}}}}^{^{\sqrt{1-{\vert x-\ov{x}\vert^2\over t^2}}}}
\ {\cos(at\ome_{2m})\over\sqrt{1-{\vert x-\ov{x}\vert^2\over t^2}-
\ome_{2m}^2}}d\ome_{2m}\bigg]\ .\cr}$$
Recall that for every $c>0$, 
$$\int_{-c}^c {\cos(at\ome)d\ome\over\sqrt{c^2-\ome^2}}=\pi J_0(atc)$$ 
where $J_0$ is the Bessel function of the first kind of order $0$. Hence
$$\eqalign{&\cos\left(\sqrt{-K}\,t\right)={1\over 2(2\pi)^{m-1}}\ppt
\ltptr^{m-1}\cr
\noalign{\vskip 5pt}
&\qquad\times\left[J_0\left(a\sqrt{t^2-\vert x-\ov{x}\vert^2}\right)
H(t^2-\vert x-\ov{x}\vert^2)\right]\ .\cr}\eqno (4.9)$$
We may continue analytically the formulas (4.8) and (4.9) in $a$ to
get (for $a$ purely imaginary) the solution of the initial value
problem for the damped wave equation ${\ptl^2 u\over\ptl t^2}=\Del u+
a^2u$, compare \cite{H}.\par
\medskip\noindent
{\bf Example~4.2}\quad {\em The harmonic oscillator\/}.\par 
Consider the operator $P=-{d^2\over dx^2}+x^2$. We wish to represent 
the solution of the initial value problem for the second order 
hyperbolic equation
$${\ptl^2 u\over\ptl t^2}=-{\ptl^2 u\over\ptl x^2}+x^2u\eqno (4.10)$$
and for this purpose we express the distribution kernel of $\cos\left(
\sqrt{P}\,t\right)=$\linebreak
$\cos\left(\sqrt{-{d^2\over dx^2}+x^2}\,t\right)$ using the formula 
(3.4) with $A={1\over i}{d\over dx}$, $B=x$. Note that $\cos(B\rho)(x,
\ov{x})=\cos(x\rho)\del(x-\ov{x})$. Hence (formally) 
$$\eqalign{&\cos\left(\sqrt{-{d^2\over dx^2}+x^2}\,t\right)=\lim_{m\rar
\infty}{1\over(2\pi)^m}\ppt\ltptr^{m-1}\bigg[t^{2m-1}\int\cdots\int\cr
\noalign{\vskip 5pt}
&\qquad\times\int_{\ome_1^2+\cdots+\ome_{2m}^2\le 1}{\del\left(x_1
-\ov{x}-{\ome_1 t\over\sqrt{m}}\right)\cos\left({\ome_2t\over\sqrt{m}}
x_2\right)\del(x_2-x_1)\over\sqrt{1-(\ome_1^2+\cdots+\ome_{2m}^2)}}\cr
\noalign{\vskip 5pt}
&\qquad\times\del\left(x_3-x_2-{\ome_3 t\over\sqrt{m}}\right)\cos\left(
{\ome_4 t\over\sqrt{m}}x_4\right)\del(x_4-x_3)\cr
\noalign{\vskip 5pt}
&\qquad\times\del\left(x_{2m-1}-x_{2m-2}-{\ome_{2m-1}\over\sqrt{m}}
\right)\cos\left({\ome_{2m} t\over\sqrt{m}}x\right)\cr
\noalign{\vskip 5pt}
&\qquad\times\del(x-x_{2m-1})d\ome_1\cdots d\ome_{2m} dx_1\cdots
dx_{2m-1}\bigg]\ .\cr}$$
Hence
$$\eqalign{&\cos\left(\sqrt{-{d^2\over dx^2}+x^2+t}\,\right)=\lim_{m\rar
\infty}{1\over(2\pi)^m}\ppt\ltptr^{m-1}\bigg[t^{2m-1}\cr
\noalign{\vskip 5pt}
&\qquad\times\int_{\ome_1^2+\cdots+\ome_{2m}^2\le 1}\cos\bigg[{\ome_2
\over\sqrt{m}}\left(\ov{x}+{\ome_1 t\over\sqrt{m}}\right)\bigg]\cos
\bigg[{\ome_2 t\over\sqrt{m}}\left(\ov{x}+{(\ome_1+\ome_3)t\over\sqrt{m}}
\right)\bigg]\cr
\noalign{\vskip 5pt}
&\qquad\times\cdots\cos\left({\ome_{2m} t\over\sqrt{m}}x\right)\del
\left(x-{\ome_1+\ome_3+\cdots+\ome_{2m-1}\over\sqrt{m}}t-\ov{x}\right)
d\ome_1\cdots d\ome_{2m}\bigg]\ .\cr}\eqno (4.11)$$
We may interpret (4.11) according to Proposition~3.1, as $P$ is
essentially self-adjoint in $L^2(R^1)$, (${1\over i}{d\over dx}$ and
$x$ are self-adjoint in their respective domains). Alternatively, set
$H=L^2(R^1;e^{-x^2}dx)$ and let $H_n(x)$ denote the Hermite polynomial 
$(-1)^ne^{x^2}{d^n\over dx^n}e^{-x^2}$. Then the well-known formulas
$H_n^{\prime}(x)=2nH_{n-1}(x)$, $xH_n(x)={H_{n+1}(x)\over 2}+nH_{n-1}
(x)$ and $\int H_n^2(x)e^{-x^2}dx=\sqrt{\pi}\,2^nn!$ imply that if 
$h\in\hbox{Sp}[H_0,H_1,\ldots,H_n]$, then $h\in D\left(A^{\alp_1}
B^{\alp_2}\cdots A^{\alp_{2m-1}}B^{\alp_{2m}}\right)$ for all $\alp$ 
and there exist constants $C$ and $K$ such that (3.8) holds. Hence 
Theorem~3.2 is applicable.\par
\medskip\noindent
{\bf Example 4.3}\quad {\em Sum of squares of vector fields\/}.\par 
Let $X$, $Y$ be vector fields with analytic coefficients defined on an
open subset $\Ome$ of $R^n$ such that $X^{\ast}=-X$, $Y^{\ast}=-Y$
(the case of more than two fields may be treated similarly, see
Remark~3.4). Let $\varphi(t,\ov{x})$, $\psi(t,\ov{x})$ denote the
solutions of the ode's
$$\eqalign{{d\varphi\over dt}(t,\ov{x})&=X\varphi(t,\ov{x}),\ \varphi(0,
\ov{x})=\ov{x}\cr
\noalign{\vskip 5pt}
{d\psi\over dt}(t,\ov{x})&=Y\psi(t,\ov{x}),\ \psi(0,\ov{x})=\ov{x}\cr}
\eqno (4.12)$$
defined in a neighborhood of $\{0\}\times\Ome$. We have, at least
formally, the relation
$$\eqalign{&\cos\left(\sqrt{X^2+Y^2}\,t\right)(x,\ov{x})=\lim_{m\rar
\infty}{1\over(2\pi)^m}\ppt\ltptr^{m-1}\bigg[t^{2m-1}\int_{\ome_1^2+
\cdots+\ome_{2m}^2}\cr
\noalign{\vskip 5pt}
&\times{\del\Big(x-\psi\Big({\ome_{2m}t\over\sqrt{m}},\varphi\Big(
{\ome_{2m-1}t\over\sqrt{m}},\psi\Big({\ome_{2m-2} t\over\sqrt{m}},
\ldots,\psi\Big({\ome_2 t\over\sqrt{m}},\varphi\Big({\ome_1 t\over
\sqrt{m}}, \ov{x}\Big)\cdots\Big)\over\sqrt{1-(\ome_1^2+\cdots+
\ome_{2m}^2)}}\cr 
\noalign{\vskip 5pt}
&\times d\ome_1\cdots d\ome_{2m}\bigg]\ .\cr}\eqno (4.13)$$
Under further assumptions the integral in the right hand side of (4.13)
will be defined for all $t$. A possible interpretation of (4.13) via
Theorem~3.2 is possible along the lines of \cite{N}.\par
Consider the case $X={\ptl\over\ptl x_1}$, $Y=x_1{\ptl\over\ptl x_2}$.
Then $X^2+Y^2$ is a self-adjoint hypoelliptic operator. Moreover, 
$\varphi(t,\ov{x})=(\ov{x}_1+t,\ov{x}_2)$, $\psi(t,\ov{x})=(\ov{x}_1,
\ov{x}_2+t)$, and we get from (4.13) that
$$\eqalign{&\cos\left(\sqrt{{\ptl^2\over\ptl x_1^2}+x_1^2{\ptl^2\over
\ptl x_2^2}}\,t\right)(x,\ov{x})=\lim_{m\rar\infty}{1\over(2\pi)^m}\ppt
\ltptr^{m-1}\bigg[t^{2m-1}\int_{\ome_1^2+\cdots+\ome_{2m}^2}\cr
\noalign{\vskip 5pt}
&\times {\del\bigg(x_1-\ov{x}_1-{t\over\sqrt{m}}\bigg(\ds{\sum_{i=1}^m}
\ome_{2i-1}\bigg),x_2-\ov{x}_2-{t\over\sqrt{m}}\bigg(\ds{\sum_{i=1}^m}
\ome_{2i}\bigg)\ov{x}_1-{t^2\over m}\sum_{1\le i\le j\le m}\ome_{2i-1}
\ome_{2j}\bigg)\over\sqrt{1-(\ome_1^2+\cdots+\ome_{2m}^2)}}\cr 
\noalign{\vskip 5pt}
&\times d\ome_1\cdots d\ome_{2m}\bigg]\ .\cr}\eqno (4.14)$$
Similarly, set $(\hbox{in}\ R^3)$ $X={\ptl\over\ptl x_1}+2x_2{\ptl\over
\ptl x_3}$, $Y={\ptl\over\ptl x_2}-2x_1{\ptl\over\ptl x_3}$. Then 
$\Del_H=X^2+Y^2$ is the Laplacian on the Heisenberg group. Here
$$\varphi=(\ov{x}_1+t,\ov{x}_2,\ov{x}_3+2\ov{x}_2t),\qquad\psi=(
\ov{x}_1,\ov{x}_2+t,\ov{x}_3-2\ov{x}_1t)\ .$$
The operators $X$, $Y$ and $\Del_H$ are left invariant. Hence it
suffices to compute kernels for $\ov{x}=0$.\par
After an elementary calculation, we get from (4.13) that
$$\eqalign{&\cos\left(\sqrt{-\Del_H}\,t\right)(x,0)=\lim_{m\rar\infty}
{1\over(2\pi)^m}\ppt\ltptr^{m-1}\bigg[t^{2m-1}\int_{\ome_1^2+\cdots+
\ome_{2m}^2}\cr
\noalign{\vskip 5pt}
&\times {\del\bigg(x_1-{t\over\sqrt{m}}\ds{\sum_{i=1}^m}\ome_{2i-1},x_2-
{t\over\sqrt{m}}\ds{\sum_{i=1}^m}\ome_{2i},x_3-{2t^2\over\sqrt{m}}\bigg(
\ds{\sum_{1\le i\le j\le m-1}}\ome_{2i}\ome_{2j-1}-\ds{\sum_{1\le i\le
j\le m}}\ome_{2i-1}\ome_{2j}\bigg)\bigg)\over\sqrt{1-(\ome_1^2+\cdots+
\ome_{2m}^2)}}\cr
\noalign{\vskip 5pt}
&\times d\ome_1\cdots d\ome_{2m}\bigg]\ .\cr}\eqno (4.15)$$
The formulas (4.13)--(4.15) may be regarded as ``Feynman-Kac''
formulas. An explicit evaluation of (4.14) and (4.15) appears to be a
challenge.\par


\begin{thebibliography}{RS}
\bibitem{A} Amoto, K., On a unitary version of Suzuki's exponential
product formula, {\em J.\ Math.\ Soc.\ Japan\/} {\bf 48} (1996)
493--499. 
\bibitem{F} Folland, G.B., {\em Introduction to Partial Differential 
Equations\/}, 2nd Edition, Princeton University Press, Princeton, NJ,
1995. 
\bibitem{H} Hadamard, J., {\em Lectures on Cauchy's Problem in Linear
Partial Differential Equations\/}, Dover Publications, New York, NY, 
1953.
\bibitem{K} Kannai, Y., Off diagonal short time asymptotics for
fundamental solutions of diffusion equations, {\em Comm.\ Part.\ Diff.\
Eq.}\ {\bf 2} (1977) 781--830.
\bibitem{N} Nelson, E., Analytic vectors, {\em Ann.\ Math.}\ {\bf 70}
(1959) 572--615.
\bibitem{RS} Reed, M.\ and Simon, B., {\em Methods of Modern 
Mathematical Physics. I.~Functional Analysis\/}, Reviewed and Enlarged
Edition, Academic Press, New York-London, 1980.
\bibitem{S} Suzuki, M., Convergence of exponential product formulas
for unbounded operators, {\em Rev.\ Math.\ Phys.}\ {\bf 8} (1996)
487--502. 
\bibitem{W} Wilks, S., {\em Mathematical Statistics\/}, Wiley, New
York-London, 1962.
\end{thebibliography}
\end{document}